\documentclass[12pt]{article}

\def\pf{\noindent {\it Proof.} }
\newcommand{\qed}{{\hfill\rule{4pt}{7pt}}}

\makeatletter
\newfont{\footsc}{cmcsc10 at 8truept}
\newfont{\footbf}{cmbx10 at 8truept}
\newfont{\footrm}{cmr10 at 10truept}
\makeatother 

\usepackage{amsmath}
\usepackage{amssymb}
\usepackage{amsfonts}
\usepackage{latexsym}

\newtheorem{thm}{Theorem}[section]
\newtheorem{rmk}{Remark}[section]
\newtheorem{prop}[thm]{Proposition}

\newtheorem{cor}[thm]{Corollary}
\newtheorem{defn}[thm]{Definition}
\newtheorem{lem}[thm]{Lemma}

\numberwithin{equation}{section}

\def\cro{\rm{cr}}
\def\ne{\rm{ne}}
\def\al{\rm{al}}
\def\sg{\rm{sg}}

\def\bl{\rm{bl}}
\def\tran{\rm{tr}}
\def\ed{\rm{ed}}

\def\M{\mathcal{M}}
\def\O{\mathcal{O}}
\def\C{\mathcal{C}}
\def\S{\mathcal{S}}
\def\T{\mathcal{T}}
\begin{document}
\begin{center}
{\Large\bf Distribution of crossings, nestings and alignments of two edges
          in matchings and partitions}
\end{center}

\vskip 2mm
\centerline{Anisse Kasraoui and Jiang Zeng}
\begin{center}
Institut Camille Jordan,
Universit\'e Claude Bernard (Lyon I)\\
F-69622, Villeurbanne Cedex, France \\
{\tt anisse@math.univ-lyon1.fr, zeng@math.univ-lyon1.fr}
\end{center}

\vskip 0.7cm \noindent{\bf Abstract.}
We construct  an involution on set partitions which keeps track of
the numbers of crossings, nestings and alignments of two edges.
 We derive then the symmetric distribution of the numbers of  crossings and nestings
 in partitions, which generalizes  Klazar's  recent result in perfect matchings.
By factorizing our involution through bijections between set partitions and some path diagrams
we obtain
the continued fraction expansions of  the corresponding ordinary generating
 functions.
%\vskip 0.2cm
%\noindent{\it Keywords}:
%\vskip 0.5cm \noindent{\bf MR Subject Classifications}: Primary
%05A18; Secondary 05A15.

\section{Introduction}

   A {\it partition} of $[n]$:=$\{1,2,\cdots,n\}$ is a collection
of disjoint nonempty subsets of $[n]$, called {\it blocks}, whose
union is $[n]$. A (perfect) \emph{matching} of $[2n]$ is a partition of $[2n]$ in $n$ two-element blocks.
Denote by $\Pi_{n}$  the set of the partitions of $[n]$
and by $\mathcal M_{2n}$ the set of the matchings of $[2n]$.
A  partition $\pi$ with $k$ blocks is written
$\pi=B_1-B_2-\cdots -B_k$, where the  blocks are ordered  in the
increasing order of their minimum elements and, within each block,
the elements are  written in the numerical order.

It is convenient to identify  a partition of $[n]$ with a \emph{partition graph}
on the vertex set $[n]$ such that there is an edge joining
$i$ and $j$ if and only if $i$ and $j$ are
\emph{consecutive elements} in a same block. We note such an edge $e$ as a
pair $(i,j)$ with $i<j$, and say that $i$ is the {\it left-hand endpoint}
of $e$ and $j$ is the {\it right-hand endpoint} of $e$.
A \emph{singleton}  is the element of a block which has only one element, so
a \emph{singleton} corresponds to  an  isolated vertex in the graph.
Conversely, a graph on the vertex set $[n]$ is a partition graph if and only if each vertex is the
left-hand (resp. right-hand) endpoint of at most one edge.
 By convention,  the vertices $1,2,\cdots,n$ are arranged on a line  in the
increasing order from left to right and an edge $(i,j)$ is drawn as an arc above the line.
An illustration is given in Figure~1.
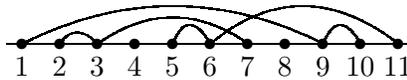
\begin{figure}[h]
\begin{center}
{\setlength{\unitlength}{1mm}
\begin{picture}(50,2.2)(0,0)
\put(-2,0){\line(1,0){54}}
\put(0,0){\circle*{1,3}}\put(0,0){\makebox(0,-6)[c]{\small 1}}
\put(5,0){\circle*{1,3}}\put(5,0){\makebox(0,-6)[c]{\small 2}}
\put(10,0){\circle*{1,3}}\put(10,0){\makebox(0,-6)[c]{\small 3}}
\put(15,0){\circle*{1,3}}\put(15,0){\makebox(0,-6)[c]{\small 4}}
\put(20,0){\circle*{1,3}}\put(20,0){\makebox(0,-6)[c]{\small 5}}
\put(25,0){\circle*{1,3}}\put(25,0){\makebox(0,-6)[c]{\small 6}}
\put(30,0){\circle*{1,3}}\put(30,0){\makebox(0,-6)[c]{\small 7}}
\put(35,0){\circle*{1,3}}\put(35,0){\makebox(0,-6)[c]{\small 8}}
\put(40,0){\circle*{1,3}}\put(40,0){\makebox(0,-6)[c]{\small 9}}
\put(45,0){\circle*{1,3}}\put(45,0){\makebox(0,-6)[c]{\small 10}}
\put(50,0){\circle*{1,3}}\put(50,0){\makebox(0,-6)[c]{\small 11}}
\qbezier(0,0)(20,10)(40,0) \qbezier(40,0)(42,5)(45,0)
\qbezier(5,0)(7,3)(10,0)\qbezier(20,0)(22,5)(25,0)
\qbezier(10,0)(20,7)(30,0)\qbezier(25,0)(37,10)(50,0)
\end{picture}
}
\end{center}
\caption{\small{Graph of the partition $\pi=\{1,9,10\}-\{2,3,7\}-\{4\}-\{5,6,11\}-\{8\}$}}
\end{figure}
\newpage
Given  a partition $\pi$ of $[n]$,
two edges $e_1=(i_{1},j_{1})$ and $e_2=(i_{2},j_{2})$
of  $\pi$ is said to form:

\begin{itemize}
\item[(i)] a $crossing$  with
$e_1$ as the $initial$ $edge$  if $i_{1}<i_{2}<j_{1}<j_{2}$ ;
\item[(ii)]  a $nesting$ with $e_2$ as
$interior$ $edge$ if $i_{1}<i_{2}<j_{2}<j_{1}$;
\item[(iii)]  an $alignment$ with
$e_1$ as $initial$  $edge$
if $i_{1}<j_{1}\leq i_{2}<j_{2}$.
\end{itemize}
An illustration of these notions is given  in Figure~2.
\begin{figure}[h]
\begin{center}
{\setlength{\unitlength}{1mm}
\begin{picture}(150,20)(-25,-10)
\put(-30,0){\line(1,0){28}}\put(10,0){\line(1,0){28}}\put(50,0){\line(1,0){28}}\put(90,0){\line(1,0){28}}
\put(-28,0){\circle*{1}}\put(-28,0){\makebox(0,-6)[c]{\small $i_1$}}
\put(-20,0){\circle*{1}}\put(-20,0){\makebox(0,-6)[c]{\small $i_2$}}
\put(-12,0){\circle*{1}}\put(-12,0){\makebox(0,-6)[c]{\small $j_1$}}
\put(-4,0){\circle*{1}}\put(-4,0){\makebox(0,-6)[c]{\small $j_2$}}
\put(12,0){\circle*{1}}\put(12,0){\makebox(0,-6)[c]{\small $i_1$}}
\put(20,0){\circle*{1}}\put(20,0){\makebox(0,-6)[c]{\small $i_2$}}
\put(28,0){\circle*{1}}\put(28,0){\makebox(0,-6)[c]{\small $j_2$}}
\put(36,0){\circle*{1}}\put(36,0){\makebox(0,-6)[c]{\small $j_1$}}
\put(52,0){\circle*{1}}\put(52,0){\makebox(0,-6)[c]{\small $i_1$}}
\put(60,0){\circle*{1}}\put(60,0){\makebox(0,-6)[c]{\small $j_1$}}
\put(68,0){\circle*{1}}\put(68,0){\makebox(0,-6)[c]{\small $i_2$}}
\put(76,0){\circle*{1}}\put(76,0){\makebox(0,-6)[c]{\small $j_2$}}
\put(92,0){\circle*{1}}\put(92,0){\makebox(0,-6)[c]{\small $i_1$}}
\put(104,0){\circle*{1}}\put(104,0){\makebox(0,-6)[c]{\small $j_1=i_2$}}
\put(116,0){\circle*{1}}\put(116,0){\makebox(0,-6)[c]{\small $j_2$}}
\qbezier(-28,0)(-20,12)(-12,0) \qbezier(-20,0)(-12,12)(-4,0)
\qbezier(12,0)(24,12)(36,0) \qbezier(20,0)(24,8)(28,0)
\qbezier(52,0)(56,12)(60,0) \qbezier(68,0)(72,12)(76,0)
\qbezier(92,0)(98,12)(104,0) \qbezier(104,0)(110,12)(116,0)
\put(-16,-10){\makebox(0,0) [c]{\small$(i)$}}
\put(24,-10){\makebox(0,0)[c]{\small$(ii)$}}
\put(64,-10){\makebox(0,0)[c]{\small$(iii)$}}
\put(104,-10){\makebox(0,0)[c]{\small$(iii')$}}
\end{picture}}
\end{center}
\caption{\small{Crossing, nesting and alignments of two edges}}
\end{figure}
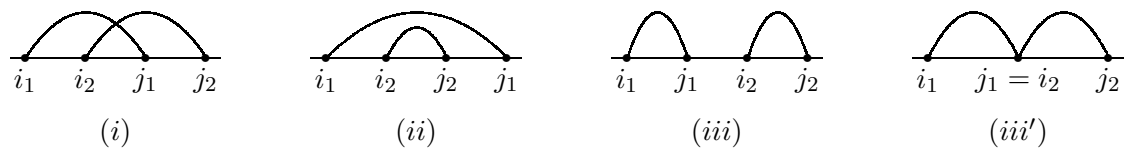\\
We denote by $\cro(\pi)$, $\ne(\pi)$ and $\al(\pi)$ the numbers of crossings,
nestings and alignments of two edges in  $\pi$, respectively.
Furthermore, consider a  block $B$ of $\pi$ whose cardinal is $\geq2$.
An element of $B$
is:
\begin{itemize}
\item[(i)] an {\it  opener}  if it is the  least element of $B$,
\item[(ii)] a {\it closer}  if it is the greatest element of $B$,
\item[(iii)]a {\it transient}  if it is neither the least nor greatest
  elements of $B$.
\end{itemize}
In the graph of $\pi$, the edges around
an opener, closer, singleton or transient are illustrated  in Figure~3.

\begin{figure}[h]
{\setlength{\unitlength}{0.80 mm}
\begin{center}
\begin{picture}(120,10)(10,0)
\put(-20,0){\line(1,0){15}}\put(-12,0){\circle*{1,3}}
\put(30,0){\line(1,0){15}}\put(38,0){\circle*{1,3}}
\put(80,0){\line(1,0){15}}\put(88,0){\circle*{1,3}}
\put(130,0){\line(1,0){15}}\put(138,0){\circle*{1,3}}
\qbezier(-12,0)(-9,5)(-5,7)\qbezier(38,0)(35,5)(31,7)
\qbezier(138,0)(135,5)(131,7)\qbezier(138,0)(141,5)(145,7)
\end{picture}
\end{center}}
\caption{\small opener, closer, singleton and transient in a partition graph}
\end{figure}
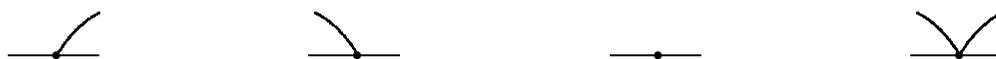
The sets of openers, closers, singletons and  transients of $\pi$
 will be denoted  by $\O (\pi)$, $\C(\pi)$, $\S(\pi)$ and $\T(\pi)$, respectively.
The $4$-tuple $\lambda(\pi)=
 (\O(\pi), \C(\pi), \S(\pi), \T(\pi))$ is called  the \emph{type} of $\pi$.

For the partition $\pi$ in Figure~1, we have $\cro(\pi)=2$, $\ne(\pi)=5$ and $\al(\pi)=8$.
Moreover, $\mathcal{O}(\pi)=\{1,2,5\}$, $\mathcal{C}(\pi)=\{7,10,11\}$, $\mathcal{S}(\pi)=\{4,8\}$
and $\mathcal{T}(\pi)=\{3,6,9\}$.
\goodbreak

\begin{defn}
A 4-tuple $\lambda=(\O,\C,\S,\T)$ of subsets of $[n]$
is a {\it partition type} of $[n]$ if there exists a partition of $[n]$
whose type is $\lambda$.
Denote by $\Pi_{n}(\lambda)$ the set
of partitions of type $\lambda$, i.e.,
$$
\Pi_{n}(\lambda)=\{\pi\in\Pi_{n}:\;\lambda(\pi)=\lambda\}.
$$
In particular,  a partition type $\lambda$ is a {\it matching type}
if $\lambda=(\O,\C):=(\O,\C,\emptyset,\emptyset)$.
Denote by $\mathcal M_{2n}(\gamma)$ the set
of matchings of type $\gamma$, i.e.,
$$
\mathcal M_{2n}(\gamma)=\{\alpha \in \mathcal M_{2n}:\;\mathcal{O}(\alpha)=\O\;and\;\mathcal{C}(\alpha)=\C\}.
$$
\end{defn}

Klazar~\cite{Kl} has recently  proved
the symmetric distribution of the numbers of crossings and nestings of two edges in perfect matchings.
The aim of this paper is to show that a much stronger result exists in the partitions which reduces to
 that of Klazar in the case of matchings.
 Note that Chen et al~\cite{Chen}
 have found other interesting results on the
crossings  and nestings  in matching and partitions, while  Corteel~\cite{Co} has given an analogous
result for permutations. Moreover we refer the reader to  Krattenthaler's
recent paper~\cite{Kr} for a more general context of related problems.

Our main result is the
 construction of an  explicit involution on the set of partitions $\Pi_n$.
\begin{thm}\label{thm:fundamental}
For each partition type $\lambda$ of $[n]$ there is an involution $\varphi:
\Pi_n(\lambda)\to \Pi_n(\lambda)$ preserving  the number of alignments,
and exchanging the numbers of crossings and nestings.
In other words, for each $\pi\in\Pi_n$, we have $\lambda(\pi)=\lambda(\varphi(\pi))$ and
\begin{align}\label{eq:property}
 \al(\varphi(\pi))=\al(\pi),\;
 \cro(\varphi(\pi))=\ne(\pi),\;\ne(\varphi(\pi))=\cro(\pi).
\end{align}
\end{thm}

\begin{cor}\label{cor:particulier}
 For each partition type $\lambda$ of $[n]$, we have
\begin{align}
 \sum_{\pi\in \Pi_{n}(\lambda)}
 p^{\cro(\pi)}q^{\ne(\pi)}t^{\al(\pi)}
 =\sum_{\pi\in\Pi_{n}(\lambda)}
 p^{\ne(\pi)}q^{\cro(\pi)}t^{\al(\pi)},\label{eq:particulier}
\end{align}
and for each matching type $\gamma$ of $[2n]$,
\begin{align}
\sum_{\alpha\in \mathcal M_{2n}(\gamma)}p^{\cro(\alpha)}q^{\ne(\alpha)}t^{\al(\alpha)}=
\sum_{\alpha\in \mathcal M_{2n}(\gamma)}p^{\ne(\alpha)}q^{\cro(\alpha)}t^{\al(\alpha)}.
 \label{eq:matching1}
 \end{align}
\end{cor}
Summing over all partition types $\lambda$ or matching types $\gamma$ we get

\begin{cor}
\begin{align}
\sum_{\pi\in \Pi_{n}}p^{cr(\pi)}q^{ne(\pi)}t^{al(\pi)}=\sum_{\pi\in
 \Pi_{n}}p^{ne(\pi)}q^{cr(\pi)}t^{al(\pi)},\label{eq:general}
\end{align}
and
\begin{align}
\sum_{\alpha\in \mathcal M_{2n}}p^{\cro(\alpha)}q^{\ne(\alpha)}t^{\al(\alpha)}=
\sum_{\alpha\in \mathcal M_{2n}}p^{\ne(\alpha)}q^{\cro(\alpha)}t^{\al(\alpha)}.\label{eq:matching2}
\end{align}
\end{cor}

In particular, by taking $t=1$ in the above corollary, we obtain
\begin{cor}\label{prop:klazar}
\begin{align}
\sum_{\pi\in \Pi_{n}}p^{cr(\pi)}q^{ne(\pi)}=\sum_{P\in
\Pi_{n}}p^{ne(\pi)}q^{cr(\pi)},\label{eq:general2}
\end{align}
and
\begin{equation}
\sum_{\alpha\in \mathcal M_{2n}}p^{\cro(\alpha)}q^{\ne(\alpha)}
=\sum_{\alpha\in \mathcal M_{2n}}p^{\ne(\alpha)}q^{\cro(\alpha)}.
\label{eq:klazar}
\end{equation}
\end{cor}
Identity~\eqref{eq:klazar} is due to Klazar~\cite{Kl}.
The $p=1$ case of \eqref{eq:klazar} had been previously  proved  by
M. de Sainte-Catherine~\cite{Sai} and by De M\'edicis and Viennot~\cite{DV}.

Our approach can be considered as an application of the
combinatorial  theory  of orthogonal polynomials developed by
Viennot~\cite{Vi} and  Flajolet~\cite{Fl}.
In fact, our involution $\varphi$ is a generalization of that
used by  De M\'edicis and Viennot~\cite{DV} for matchings.
A variant of this bijection on partitions has been used
by Ksavrelof and Zeng~\cite{KZ} to prove other equinumerous results on partitions.

The paper is organized as follows: we present the
involution $\varphi$ and the proof of theorem 1.1 in
section 2; in section 3 we factorize our involution through two bijections
 between partitions and Charlier diagrams, which permit us to derive continued fraction expansions of
the ordinary generating functions with respect to the numbers of crossings and nestings
of two edges in matchings and partitions.

\section{Proof of Theorem~\ref{thm:fundamental}}
Let $\pi=B_1-B_2-\cdots -B_k$ be a partition of $[n]$ and $i$ an integer in $[n]$.
The restriction  $B_j(\leq i):=B_j\cap [i]$ of the block $B_j$
is  said to be
\emph{opened} (resp. \emph{closed} and  \emph{empty} ) if $B\not\subseteq [i]$ (resp.
   $B\subseteq [i]$ and $B\cap [i]=\emptyset$). The $i$-th \emph{trace}  of $\pi$  is
   defined by
 $$
 T_i(\pi)=B_1(\leq i)-B_2(\leq i)-\cdots -B_k(\leq
i).
$$
We can represent $T_i(\pi)$ by a graph $D_i(\pi)$ on the vertex
set $[i]$.  Define  $D_i(\pi)$  as the subgraph of the graph of
$\pi$
 induced by the vertex set $[i]$, with
the additional condition that
 for any edge
 $(x,y)$ of $\pi$ such that $x\leq i<y$, we attach a "half-edge" to the vertex $x$,
 called \emph{vacant vertex}. Denote by $l_{i}(\pi)$ the number of vacant vertices in
$D_{i-1}(\pi)$, with $D_0=\emptyset$. Moreover, if $i$ is a closer
or a transient, there is an edge
 $(j,i)$ with  $j<i$, we denote by $\gamma_i(\pi)$  the rank of the vertex $j$
 among the vacant vertices
of $D_{i-1}(\pi)$, the vacant vertices being arranged from left to
right in the order of their creation, namely, in increasing order.

For instance, if $\pi$ is the partition given in Figure~1, then
$T_6(\pi)=\{1,\ldots \}-\{2,3,\ldots \}-\{4\}-\{5,6,\ldots\}$,
where  each opened block has an ellipsis.
  The corresponding graphs $D_5(\pi)$ and $D_6(\pi)$  are presented in
  Figure~4. We have $l_6(\pi)=3$
and $\gamma_6(\pi)=3$.
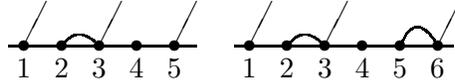
\begin{figure}[h]
\begin{center}
%%%%%%%%%%%%%%%%%%%%%%%%%%%%%%%%%%%%%%%%%%%%%%%%%%%%%%%%%%%%%%%%%%%%%%%%%%%%%%D5
{\setlength{\unitlength}{1mm}
\begin{picture}(60,5)(-5,3)
\put(-2,0){\line(1,0){25}}
\put(0,0){\circle*{1,3}}\put(0,0){\makebox(0,-6)[c]{\small 1}}
\put(5,0){\circle*{1,3}}\put(5,0){\makebox(0,-6)[c]{\small 2}}
\put(10,0){\circle*{1,3}}\put(10,0){\makebox(0,-6)[c]{\small 3}}
\put(15,0){\circle*{1,3}}\put(15,0){\makebox(0,-6)[c]{\small 4}}
\put(20,0){\circle*{1,3}}\put(20,0){\makebox(0,-6)[c]{\small 5}}
\qbezier(5,0)(7,3)(10,0)
\put(0,0){\line(1,2){3}}\put(10,0){\line(1,2){3}}\put(20,0){\line(1,2){3}}

%%%%%%%%%%%%%%%%%%%%%%%%%%%%%%%%%%%%%%%%%%%%%%%%%%%%%%%%%%%%%%%%%%%%%%%%%%%%%%%%%D6
\put(28,0){\line(1,0){29}}
\put(30,0){\circle*{1,3}}\put(30,0){\makebox(0,-6)[c]{\small 1}}
\put(35,0){\circle*{1,3}}\put(35,0){\makebox(0,-6)[c]{\small 2}}
\put(40,0){\circle*{1,3}}\put(40,0){\makebox(0,-6)[c]{\small 3}}
\put(45,0){\circle*{1,3}}\put(45,0){\makebox(0,-6)[c]{\small 4}}
\put(50,0){\circle*{1,3}}\put(50,0){\makebox(0,-6)[c]{\small 5}}
\put(55,0){\circle*{1,3}}\put(55,0){\makebox(0,-6)[c]{\small 6}}
\qbezier(35,0)(37,3)(40,0)\qbezier(50,0)(52,5)(55,0)
\put(30,0){\line(1,2){3}}\put(40,0){\line(1,2){3}}
\put(55,0){\line(1,2){3}}
\end{picture}
}
\end{center}
\caption{\small Graphs of $D_5(\pi)$ and $D_6(\pi)$}
\end{figure}

\goodbreak
Now, we can describe our fundamental bijection $\varphi$ using partition
graphs.
 In the following,
by "declare the vertex $i$  vacant" we mean "attach a half-edge to
the vertex $i$".
 Let $\pi\in \Pi_n$, with type $\lambda=(\O,\C,\S,\T)$. We obtain $\varphi(\pi)$ by the following algorithm:
\begin{enumerate}
\item Set $D'_0=\emptyset$.
\item For $1\leq i\leq n$, the graph $D'_i$ is
obtained from $D'_{i-1}$ by adding $i$ as follows:
 \begin{itemize}
\item [(i)] if $i\in O$, declare the vertex $i$  vacant.
 \item[(ii)] if $i\in S$, add $i$ as an isolated vertex.
 \item[(iii)] if $i\in C\cup T$, join $i$ to the $\gamma_i(\pi)$-th (from right to left) vacant
vertex of $D'_{i-1}$. Moreover, if $i\in T$, declare the vertex
$i$ vacant.
\end{itemize}
\item Set $\varphi(\pi):= D'_n$
\end{enumerate}

\medskip

 \begin{lem}\label{lem:a}
 The mapping $\varphi:\Pi_n\to\Pi_n$ is well defined. Moreover, it
 is an involution which preserves the type.
\end{lem}

\pf  By induction on $i$, it is easy
to see that $D'_i$ ($0\leq i\leq n$) has the same vacant vertices
as $D_i(\pi)$. So $(iii)$ is valid and  $D'_n$ is a partition
graph of $[n]$. The algorithm is well defined.
By inspecting the algorithm, we see that
$\varphi(\pi)$ has the same type as $\pi$. Finally, the operation "reverse the
order of the vacant vertices twice" preserves the
original order. So $\varphi$ is an involution.\qed

\begin{rmk}
The graph $D_i'$ corresponds with the graph of $i$-th trace of
$\varphi(\pi)$.
\end{rmk}

For instance, if
$\pi$ is that in Figure~1, then
$\varphi(\pi)=\{1,3,10\}-\{2,6,9,11\}-\{4\}-\{5,7\}-\{8\}$. Notice
that $\cro(\varphi(\pi))=\ne(\pi)=5$,
$\ne(\varphi(\pi))=\cro(\pi)=2$ and
$\al(\varphi(\pi))=\al(\pi)=8$.
 An
illustration of the step-by-step construction of $\varphi(\pi)$ is given in
Figure~5.

To complete the proof of Theorem 1.1 it remains to verify
\eqref{eq:property}. In fact we shall prove a stronger result. For
any closer or transient $j$ of a partition $\pi$, let
$\cro(\pi;j)$ (resp. $\ne(\pi;j)$ and $\al(\pi;j)$) be the number
of crossings (resp. nestings and alignments) whose initial (resp.
interior and initial) edge has
 $j$ as the right-hand endpoint. Clearly
$$
\cro(\pi)=\sum\cro(\pi;j),\qquad \ne(\pi)=\sum\ne(\pi;j),\qquad \al(\pi)=\sum\al(\pi;j),
$$
where the summations are over $j\in C(\pi)\cup T(\pi)$.

\begin{lem}\label{lem:principal}
 Let $\pi$ be a partition of $[n]$ and $j$ a closer or transient of $\pi$. Then
\begin{align*}
\al(\varphi(\pi);j)= \al(\pi;j),\quad
\cro(\varphi(\pi);j)= \ne(\pi;j),\quad
\ne(\varphi(\pi);j)= \cro(\pi;j).
\end{align*}
\end{lem}
\pf For any partition $\pi$,
the number of alignments with $j$ as the right-hand endpoint, i.e. $\al(\pi;j)$,
 is equal too the number of openers and
transients which are $\geq j$. Now, as  $\varphi(\pi)$ has the same openers and transients as
$\pi$,
we get immediately $\al(\varphi(\pi);j)=\al(\pi;j)$.

\begin{figure}[h]
\begin{center}
{\setlength{\unitlength}{0.8mm}
\begin{picture}(120,170)(-30,-30)
%%%%%%%%%%%%%%%%%%%%%%%%%%%%%%%%%%%%%%%%%%%%%%%%%%%%%%%%%%
\put(15,135){\makebox(0,0)[c]{$D_i(\pi)$}}\put(75,135){\makebox(0,0)[c]{$D'_i$}}
\put(-35,135){\makebox(0,0)[c]{$i$}}\put(-20,135){\makebox(0,0)[c]{$\gamma_i(\pi)$}}
\put(-35,-30){\makebox(0,0)[c]{\small
$11$}}\put(-20,-30){\makebox(0,0)[c]{\small $1$}}
\put(-35,-15){\makebox(0,0)[c]{\small
$10$}}\put(-20,-15){\makebox(0,0)[c]{\small $2$}}
\put(-35,0){\makebox(0,0)[c]{\small
$9$}}\put(-20,0){\makebox(0,0)[c]{\small $1$}}
\put(-35,15){\makebox(0,0)[c]{\small $8$}}
\put(-35,30){\makebox(0,0)[c]{\small
$7$}}\put(-20,30){\makebox(0,0)[c]{\small $2$}}
\put(-35,45){\makebox(0,0)[c]{\small
$6$}}\put(-20,45){\makebox(0,0)[c]{\small $3$}}
\put(-35,60){\makebox(0,0)[c]{\small $5$}}
\put(-35,75){\makebox(0,0)[c]{\small $4$}}
\put(-35,90){\makebox(0,0)[c]{\small
$3$}}\put(-20,90){\makebox(0,0)[c]{\small $2$}}
\put(-35,105){\makebox(0,0)[c]{\small $2$}}
\put(-35,120){\makebox(0,0)[c]{\small $1$}}

\put(0,-15){\circle*{1,3}}\put(0,-15){\makebox(0,-6)[c]{\small 1}}\put(60,-15){\circle*{1,3}}\put(60,-15){\makebox(0,-6)[c]{\small 1}}
\put(5,-15){\circle*{1,3}}\put(5,-15){\makebox(0,-6)[c]{\small 2}}\put(65,-15){\circle*{1,3}}\put(65,-15){\makebox(0,-6)[c]{\small 2}}
\put(10,-15){\circle*{1,3}}\put(10,-15){\makebox(0,-6)[c]{\small 3}}\put(70,-15){\circle*{1,3}}\put(70,-15){\makebox(0,-6)[c]{\small 3}}
\put(15,-15){\circle*{1,3}}\put(15,-15){\makebox(0,-6)[c]{\small 4}}\put(75,-15){\circle*{1,3}}\put(75,-15){\makebox(0,-6)[c]{\small 4}}
\put(20,-15){\circle*{1,3}}\put(20,-15){\makebox(0,-6)[c]{\small 5}}\put(80,-15){\circle*{1,3}}\put(80,-15){\makebox(0,-6)[c]{\small 5}}
\put(25,-15){\circle*{1,3}}\put(25,-15){\makebox(0,-6)[c]{\small 6}}\put(85,-15){\circle*{1,3}}\put(85,-15){\makebox(0,-6)[c]{\small 6}}
\put(30,-15){\circle*{1,3}}\put(30,-15){\makebox(0,-6)[c]{\small 7}}\put(90,-15){\circle*{1,3}}\put(90,-15){\makebox(0,-6)[c]{\small 7}}
\put(35,-15){\circle*{1,3}}\put(35,-15){\makebox(0,-6)[c]{\small 8}}\put(95,-15){\circle*{1,3}}\put(95,-15){\makebox(0,-6)[c]{\small 8}}
\put(40,-15){\circle*{1,3}}\put(40,-15){\makebox(0,-6)[c]{\small 9}}\put(100,-15){\circle*{1,3}}\put(100,-15){\makebox(0,-6)[c]{\small 9}}
\put(45,-15){\circle*{1,3}}\put(45,-15){\makebox(0,-6)[c]{\small 10}}\put(105,-15){\circle*{1,3}}\put(105,-15){\makebox(0,-6)[c]{\small 10}}

\put(0,-15){\circle*{1,3}}\put(0,-15){\makebox(0,-6)[c]{\small 1}}\put(60,-15){\circle*{1,3}}\put(60,-15){\makebox(0,-6)[c]{\small 1}}
\put(5,-15){\circle*{1,3}}\put(5,-15){\makebox(0,-6)[c]{\small 2}}\put(65,-15){\circle*{1,3}}\put(65,-15){\makebox(0,-6)[c]{\small 2}}
\put(10,-15){\circle*{1,3}}\put(10,-15){\makebox(0,-6)[c]{\small 3}}\put(70,-15){\circle*{1,3}}\put(70,-15){\makebox(0,-6)[c]{\small 3}}
\put(15,-15){\circle*{1,3}}\put(15,-15){\makebox(0,-6)[c]{\small 4}}\put(75,-15){\circle*{1,3}}\put(75,-15){\makebox(0,-6)[c]{\small 4}}
\put(20,-15){\circle*{1,3}}\put(20,-15){\makebox(0,-6)[c]{\small5}}\put(80,-15){\circle*{1,3}}\put(80,-15){\makebox(0,-6)[c]{\small5}}
\put(25,-15){\circle*{1,3}}\put(25,-15){\makebox(0,-6)[c]{\small6}}\put(85,-15){\circle*{1,3}}\put(85,-15){\makebox(0,-6)[c]{\small6}}
\put(30,-15){\circle*{1,3}}\put(30,-15){\makebox(0,-6)[c]{\small7}}\put(90,-15){\circle*{1,3}}\put(90,-15){\makebox(0,-6)[c]{\small7}}
\put(35,-15){\circle*{1,3}}\put(35,-15){\makebox(0,-6)[c]{\small8}}\put(95,-15){\circle*{1,3}}\put(95,-15){\makebox(0,-6)[c]{\small8}}
\put(40,-15){\circle*{1,3}}\put(40,-15){\makebox(0,-6)[c]{\small9}}\put(100,-15){\circle*{1,3}}\put(100,-15){\makebox(0,-6)[c]{\small9}}

\put(0,0){\circle*{1,3}}\put(0,0){\makebox(0,-6)[c]{\small 1}}\put(60,0){\circle*{1,3}}\put(60,0){\makebox(0,-6)[c]{\small1}} \put(5,30){\circle*{1,3}}\put(5,30){\makebox(0,-6)[c]{\small2}}\put(65,30){\circle*{1,3}}\put(65,30){\makebox(0,-6)[c]{\small2}}
\put(5,0){\circle*{1,3}}\put(5,0){\makebox(0,-6)[c]{\small2}}\put(65,0){\circle*{1,3}}\put(65,0){\makebox(0,-6)[c]{\small2}}
\put(10,0){\circle*{1,3}}\put(10,0){\makebox(0,-6)[c]{\small3}}\put(70,0){\circle*{1,3}}\put(70,0){\makebox(0,-6)[c]{\small3}}
\put(15,0){\circle*{1,3}}\put(15,0){\makebox(0,-6)[c]{\small4}}\put(75,0){\circle*{1,3}}\put(75,0){\makebox(0,-6)[c]{\small4}}
\put(20,0){\circle*{1,3}}\put(20,0){\makebox(0,-6)[c]{\small5}}\put(80,0){\circle*{1,3}}\put(80,0){\makebox(0,-6)[c]{\small5}}
\put(25,0){\circle*{1,3}}\put(25,0){\makebox(0,-6)[c]{\small6}}\put(85,0){\circle*{1,3}}\put(85,0){\makebox(0,-6)[c]{\small6}}
\put(30,0){\circle*{1,3}}\put(30,0){\makebox(0,-6)[c]{\small7}}\put(90,0){\circle*{1,3}}\put(90,0){\makebox(0,-6)[c]{\small7}}
\put(35,0){\circle*{1,3}}\put(35,0){\makebox(0,-6)[c]{\small8}}\put(95,0){\circle*{1,3}}\put(95,0){\makebox(0,-6)[c]{\small8}}
\put(40,0){\circle*{1,3}}\put(40,0){\makebox(0,-6)[c]{\small9}}\put(100,0){\circle*{1,3}}\put(100,0){\makebox(0,-6)[c]{\small9}}

\put(0,15){\circle*{1,3}}\put(0,15){\makebox(0,-6)[c]{\small 1}}\put(60,15){\circle*{1,3}}\put(60,15){\makebox(0,-6)[c]{\small1}} \put(5,30){\circle*{1,3}}\put(5,30){\makebox(0,-6)[c]{\small2}}\put(65,30){\circle*{1,3}}\put(65,30){\makebox(0,-6)[c]{\small2}}
\put(5,15){\circle*{1,3}}\put(5,15){\makebox(0,-6)[c]{\small2}}\put(65,15){\circle*{1,3}}\put(65,15){\makebox(0,-6)[c]{\small2}}
\put(10,15){\circle*{1,3}}\put(10,15){\makebox(0,-6)[c]{\small3}}\put(70,15){\circle*{1,3}}\put(70,15){\makebox(0,-6)[c]{\small3}}
\put(15,15){\circle*{1,3}}\put(15,15){\makebox(0,-6)[c]{\small4}}\put(75,15){\circle*{1,3}}\put(75,15){\makebox(0,-6)[c]{\small4}}
\put(20,15){\circle*{1,3}}\put(20,15){\makebox(0,-6)[c]{\small5}}\put(80,15){\circle*{1,3}}\put(80,15){\makebox(0,-6)[c]{\small5}}
\put(25,15){\circle*{1,3}}\put(25,15){\makebox(0,-6)[c]{\small6}}\put(85,15){\circle*{1,3}}\put(85,15){\makebox(0,-6)[c]{\small6}}
\put(30,15){\circle*{1,3}}\put(30,15){\makebox(0,-6)[c]{\small7}}\put(90,15){\circle*{1,3}}\put(90,15){\makebox(0,-6)[c]{\small7}}
\put(35,15){\circle*{1,3}}\put(35,15){\makebox(0,-6)[c]{\small8}}\put(95,15){\circle*{1,3}}\put(95,15){\makebox(0,-6)[c]{\small8}}

\put(0,30){\circle*{1,3}}\put(0,30){\makebox(0,-6)[c]{\small 1}}\put(60,30){\circle*{1,3}}\put(60,30){\makebox(0,-6)[c]{\small1}} \put(5,45){\circle*{1,3}}\put(5,45){\makebox(0,-6)[c]{\small2}}\put(65,45){\circle*{1,3}}\put(65,45){\makebox(0,-6)[c]{\small2}}
\put(5,30){\circle*{1,3}}\put(5,30){\makebox(0,-6)[c]{\small2}}\put(65,30){\circle*{1,3}}\put(65,30){\makebox(0,-6)[c]{\small2}}
\put(10,30){\circle*{1,3}}\put(10,30){\makebox(0,-6)[c]{\small3}}\put(70,30){\circle*{1,3}}\put(70,30){\makebox(0,-6)[c]{\small3}}
\put(15,30){\circle*{1,3}}\put(15,30){\makebox(0,-6)[c]{\small4}}\put(75,30){\circle*{1,3}}\put(75,30){\makebox(0,-6)[c]{\small4}}
\put(20,30){\circle*{1,3}}\put(20,30){\makebox(0,-6)[c]{\small5}}\put(80,30){\circle*{1,3}}\put(80,30){\makebox(0,-6)[c]{\small5}}
\put(25,30){\circle*{1,3}}\put(25,30){\makebox(0,-6)[c]{\small6}}\put(85,30){\circle*{1,3}}\put(85,30){\makebox(0,-6)[c]{\small6}}
\put(30,30){\circle*{1,3}}\put(30,30){\makebox(0,-6)[c]{\small7}}\put(90,30){\circle*{1,3}}\put(90,30){\makebox(0,-6)[c]{\small7}}

\put(0,45){\circle*{1,3}}\put(0,45){\makebox(0,-6)[c]{\small 1}}\put(60,45){\circle*{1,3}}\put(60,45){\makebox(0,-6)[c]{\small1}} \put(5,60){\circle*{1,3}}\put(5,60){\makebox(0,-6)[c]{\small2}}\put(65,60){\circle*{1,3}}\put(65,60){\makebox(0,-6)[c]{\small2}}
\put(5,45){\circle*{1,3}}\put(5,45){\makebox(0,-6)[c]{\small2}}\put(65,45){\circle*{1,3}}\put(65,45){\makebox(0,-6)[c]{\small2}}
\put(10,45){\circle*{1,3}}\put(10,45){\makebox(0,-6)[c]{\small3}}\put(70,45){\circle*{1,3}}\put(70,45){\makebox(0,-6)[c]{\small3}}
\put(15,45){\circle*{1,3}}\put(15,45){\makebox(0,-6)[c]{\small4}}\put(75,45){\circle*{1,3}}\put(75,45){\makebox(0,-6)[c]{\small4}}
\put(20,45){\circle*{1,3}}\put(20,45){\makebox(0,-6)[c]{\small5}}\put(80,45){\circle*{1,3}}\put(80,45){\makebox(0,-6)[c]{\small5}}
\put(25,45){\circle*{1,3}}\put(25,45){\makebox(0,-6)[c]{\small6}}\put(85,45){\circle*{1,3}}\put(85,45){\makebox(0,-6)[c]{\small6}}

\put(0,60){\circle*{1,3}}\put(0,60){\makebox(0,-6)[c]{\small 1}}\put(60,60){\circle*{1,3}}\put(60,60){\makebox(0,-6)[c]{\small1}} \put(5,75){\circle*{1,3}}\put(5,75){\makebox(0,-6)[c]{\small2}}\put(65,75){\circle*{1,3}}\put(65,75){\makebox(0,-6)[c]{\small2}}
\put(5,60){\circle*{1,3}}\put(5,60){\makebox(0,-6)[c]{\small2}}\put(65,60){\circle*{1,3}}\put(65,60){\makebox(0,-6)[c]{\small2}}
\put(10,60){\circle*{1,3}}\put(10,60){\makebox(0,-6)[c]{\small3}}\put(70,60){\circle*{1,3}}\put(70,60){\makebox(0,-6)[c]{\small3}}
\put(15,60){\circle*{1,3}}\put(15,60){\makebox(0,-6)[c]{\small4}}\put(75,60){\circle*{1,3}}\put(75,60){\makebox(0,-6)[c]{\small4}}
\put(20,60){\circle*{1,3}}\put(20,60){\makebox(0,-6)[c]{\small5}}\put(80,60){\circle*{1,3}}\put(80,60){\makebox(0,-6)[c]{\small5}}

\put(0,75){\circle*{1,3}}\put(0,75){\makebox(0,-6)[c]{\small 1}}\put(60,75){\circle*{1,3}}\put(60,75){\makebox(0,-6)[c]{\small1}} \put(5,90){\circle*{1,3}}\put(5,90){\makebox(0,-6)[c]{\small2}}\put(65,90){\circle*{1,3}}\put(65,90){\makebox(0,-6)[c]{\small2}}
\put(10,75){\circle*{1,3}}\put(10,75){\makebox(0,-6)[c]{\small3}}\put(70,75){\circle*{1,3}}\put(70,75){\makebox(0,-6)[c]{\small3}}
\put(15,75){\circle*{1,3}}\put(15,75){\makebox(0,-6)[c]{\small4}}\put(75,75){\circle*{1,3}}\put(75,75){\makebox(0,-6)[c]{\small4}}

\put(0,90){\circle*{1,3}}\put(0,90){\makebox(0,-6)[c]{\small1}}\put(60,90){\circle*{1,3}}\put(60,90){\makebox(0,-6)[c]{\small1}}
\put(5,90){\circle*{1,3}}\put(5,90){\makebox(0,-6)[c]{\small2}}\put(65,90){\circle*{1,3}}\put(65,90){\makebox(0,-6)[c]{\small2}}
\put(10,90){\circle*{1,3}}\put(10,90){\makebox(0,-6)[c]{\small3}}\put(70,90){\circle*{1,3}}\put(70,90){\makebox(0,-6)[c]{\small3}}

\put(0,105){\circle*{1,3}}\put(0,105){\makebox(0,-6)[c]{\small1}}\put(60,105){\circle*{1,3}}\put(60,105){\makebox(0,-6)[c]{\small1}}
\put(5,105){\circle*{1,3}}\put(5,105){\makebox(0,-6)[c]{\small2}}\put(65,105){\circle*{1,3}}\put(65,105){\makebox(0,-6)[c]{\small2}}

\put(0,120){\circle*{1,3}}\put(0,120){\makebox(0,-6)[c]{\small1}}\put(60,120){\circle*{1,3}}\put(60,120){\makebox(0,-6)[c]{\small1}}
%%%%%%%%%%%%%%%%%%%%%%%%%%%%%%%%

\put(0,15){\line (1,2){3}} \put(0,30){\line (1,2){3}}\put(0,45){\line (1,2){3}}\put(0,60){\line(1,2){3}}
\put(0,75){\line(1,2){3}} \put(0,90){\line (1,2){3}}\put(0,105){\line (1,2){3}} \put(0,120){\line (1,2){3}}

\put(5,105){\line (1,2){3}} \put(10,90){\line (1,2){3}}\put(10,75){\line (1,2){3}} \put(20,60){\line (1,2){3}}
\put(10,45){\line (1,2){3}} \put(40,0){\line(1,2){3}}\put(10,60){\line (1,2){3}}

\qbezier(5,90)(7,100)(10,90) \qbezier(5,75)(7,85)(10,75)\qbezier(5,60)(7,70)(10,60) \qbezier(5,45)(7,55)(10,45)
\qbezier(5,30)(7,40)(10,30) \qbezier(5,15)(7,25)(10,15)\qbezier(5,0)(7,7)(10,0) \qbezier(5,-15)(7,-8)(10,-15)

\qbezier(20,45)(22,55)(25,45) \qbezier(20,30)(22,40)(25,30)\qbezier(20,15)(22,25)(25,15) \qbezier(20,0)(22,10)(25,0)
\qbezier(20,-15)(22,-5)(25,-15)

\qbezier(10,30)(20,43)(30,30) \qbezier(10,15)(20,28)(30,15)\qbezier(10,0)(20,13)(30,0) \qbezier(10,-15)(20,-2)(30,-15)

\qbezier(0,0)(20,15)(40,0) \qbezier(0,-15)(20,0)(40,-15)\qbezier(40,-15)(42,-5)(45,-15)
%%%%%%%%%%%%%%%%%%%%%%%%%%%%%%%%%%%%%%%%%%%%%%%%%%%%%%

\put(25,-15){\line (1,2){4}}\put(25,0){\line (1,2){4}}\put(25,15){\line (1,2){3}}
\put(25,30){\line (1,2){3}}\put(25,45){\line (1,2){3}}

\qbezier(0,-30)(20,-17)(40,-30)\qbezier(5,-30)(8,-25)(10,-30)
\qbezier(10,-30)(20,-21)(30,-30)\qbezier(20,-30)(22,-25)(25,-30)
\qbezier(25,-30)(37,-17)(50,-30)\qbezier(40,-30)(42,-20)(45,-30)

\qbezier(60,90)(65,100)(70,90) \qbezier(60,75)(65,85)(70,75)\qbezier(60,60)(65,70)(70,60)
\qbezier(60,45)(65,55)(70,45)\qbezier(60,30)(65,40)(70,30) \qbezier(60,15)(65,25)(70,15)
\qbezier(60,0)(65,7)(70,0)\qbezier(60,-15)(65,-8)(70,-15)\qbezier(60,-30)(65,-25)(70,-30)
\qbezier(65,45)(75,55)(85,45) \qbezier(65,30)(75,40)(85,30)\qbezier(65,15)(75,25)(85,15)
\qbezier(65,0)(75,7)(85,0)\qbezier(65,-15)(75,-8)(85,-15) \qbezier(65,-30)(75,-23)(85,-30)
\qbezier(80,30)(85,35)(90,30)\qbezier(80,15)(85,20)(90,15)\qbezier(80,0)(85,5)(90,0)
\qbezier(80,-15)(85,-10)(90,-15)\qbezier(80,-30)(85,-23)(90,-30)
\qbezier(85,0)(92,10)(100,0)\qbezier(85,-15)(92,-5)(100,-15)\qbezier(85,-30)(92,-23)(100,-30)
\qbezier(70,-15)(85,-1)(105,-15)\qbezier(70,-30)(85,-17)(105,-30)
\qbezier(100,-30)(105,-16)(110,-30)

\put(0,-30){\circle*{1,3}}\put(0,-30){\makebox(0,-6)[c]{\small 1}}\put(60,-30){\circle*{1,3}}\put(60,-30){\makebox(0,-6)[c]{\small 1}}
\put(5,-30){\circle*{1,3}}\put(5,-30){\makebox(0,-6)[c]{\small 2}}\put(65,-30){\circle*{1,3}}\put(65,-30){\makebox(0,-6)[c]{\small 2}}
\put(10,-30){\circle*{1,3}}\put(10,-30){\makebox(0,-6)[c]{\small 3}}\put(70,-30){\circle*{1,3}}\put(70,-30){\makebox(0,-6)[c]{\small 3}}
\put(15,-30){\circle*{1,3}}\put(15,-30){\makebox(0,-6)[c]{\small 4}}\put(75,-30){\circle*{1,3}}\put(75,-30){\makebox(0,-6)[c]{\small 4}}
\put(20,-30){\circle*{1,3}}\put(20,-30){\makebox(0,-6)[c]{\small 5}}\put(80,-30){\circle*{1,3}}\put(80,-30){\makebox(0,-6)[c]{\small 5}}
\put(25,-30){\circle*{1,3}}\put(25,-30){\makebox(0,-6)[c]{\small 6}}\put(85,-30){\circle*{1,3}}\put(85,-30){\makebox(0,-6)[c]{\small 6}}
\put(30,-30){\circle*{1,3}}\put(30,-30){\makebox(0,-6)[c]{\small 7}}\put(90,-30){\circle*{1,3}}\put(90,-30){\makebox(0,-6)[c]{\small 7}}
\put(35,-30){\circle*{1,3}}\put(35,-30){\makebox(0,-6)[c]{\small 8}}\put(95,-30){\circle*{1,3}}\put(95,-30){\makebox(0,-6)[c]{\small 8}}
\put(40,-30){\circle*{1,3}}\put(40,-30){\makebox(0,-6)[c]{\small 9}}\put(100,-30){\circle*{1,3}}\put(100,-30){\makebox(0,-6)[c]{\small9}}
\put(45,-30){\circle*{1,3}}\put(45,-30){\makebox(0,-6)[c]{\small10}}\put(105,-30){\circle*{1,3}}\put(105,-30){\makebox(0,-6)[c]{\small10}}
\put(50,-30){\circle*{1,3}}\put(50,-30){\makebox(0,-6)[c]{\small11}}\put(110,-30){\circle*{1,3}}\put(110,-30){\makebox(0,-6)[c]{\small11}}

\put(60,120){\line (1,2){4}} \put(60,105){\line (1,2){4}}\put(65,60){\line(1,2){3}} \put(65,75){\line(1,2){3}}
\put(65,90){\line (1,2){3}} \put(65,105){\line (1,2){3}}\put(70,0){\line(1,2){3}}\put(70,15){\line (1,2){3}} \put(70,30){\line(1,2){3}}
\put(70,45){\line (1,2){3}} \put(70,60){\line(1,2){3}}\put(70,75){\line(1,2){3}} \put(70,90){\line (1,2){3}}
\put(80,45){\line (1,2){3}}\put(80,60){\line(1,2){3}}\put(100,0){\line (1,2){3}}\put(100,-15){\line (1,2){3}}
\put(85,15){\line (1,2){3}}\put(85,30){\line (1,2){3}}\put(85,45){\line (1,2){3}}
\end{picture}
}
\medskip
\end{center}
\caption {Construction of $\varphi(\pi)=\{1,3,10\}-\{2,6,9,11\}-\{4\}-\{5,7\}-\{8\}$}
\end{figure}
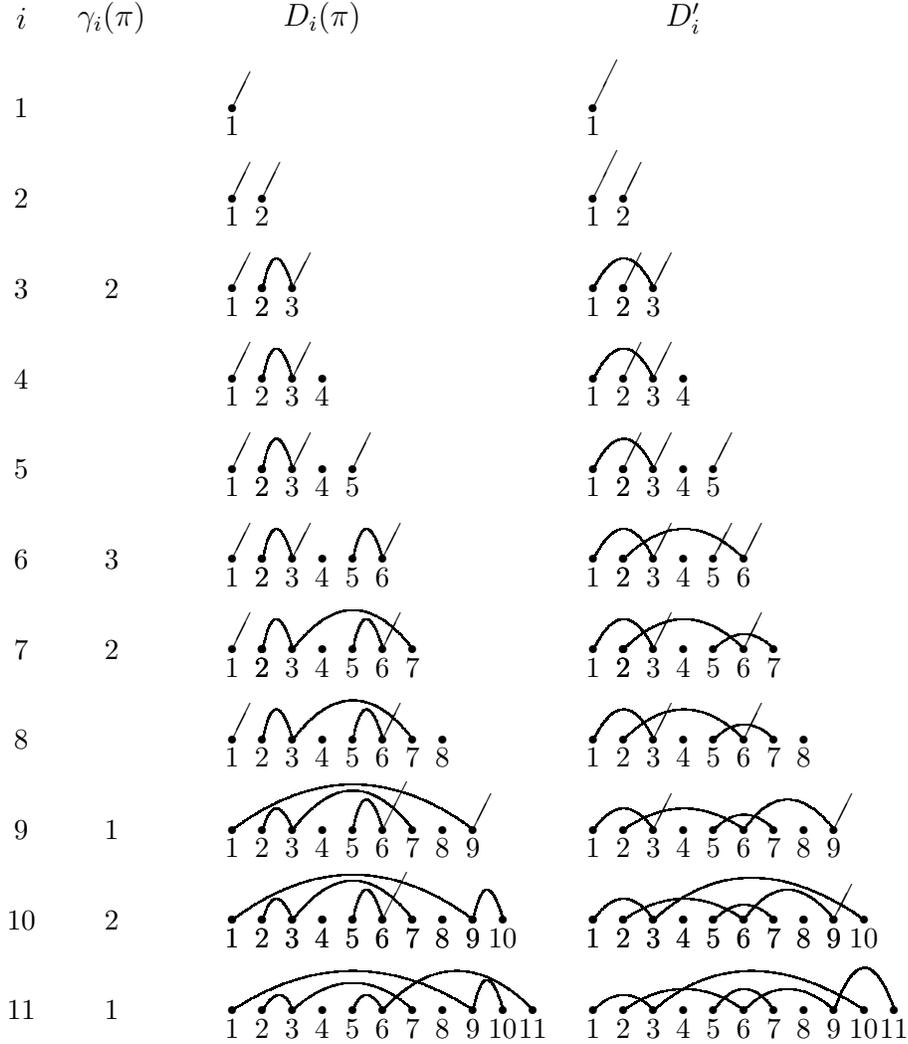

 Next, in the $j$-th ($1\leq j\leq n-1$) step of the construction of $\varphi(\pi)$, we
add the vertex $j$ to $D'_{j-1}$ for
obtaining $D'_{j}$. There are exactly $l_j:=l_j(\pi)$ vacant vertices in $D'_{j-1}$ (resp. $D_{j-1})$. These vertices are smaller than $j$
and  arranged from left to right in increasing order.
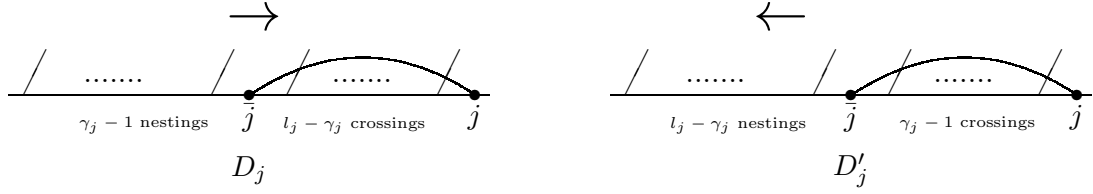
\begin{figure}[h]
\begin{center}
{\setlength{\unitlength}{1mm}
\begin{picture}(124,20)(0,-8)

\put(-2,0){\line(1,0){64}}
\put(27,10){\makebox(0,0)[l]{\LARGE{$\rightarrow$}}}
\put(0,0){\line(1,2){3}}
\put(12,2){\makebox(0,0)[c]{.......}}
\put(25,0){\line(1,2){3}}
\put(30,0){\circle*{1,3}}\put(30,0){\makebox(0,-6)[c]{\small $\bar j$}}
\put(35,0){\line(1,2){3}}
\put(45,2){\makebox(0,0)[c]{.......}}
\put(55,0){\line(1,2){3}}
\put(60,0){\circle*{1,3}}\put(60,0){\makebox(0,-6)[c]{\small $j$}}
\qbezier(30,0)(45,10)(60,0)
\put(16,-4){\makebox(0,0)[c]{\tiny $\gamma_j-1$ nestings}}
\put(43,-4){\makebox(2,0)[c]{\tiny $l_j-\gamma_j$ crossings}}

\put(78,0){\line(1,0){64}}
\put(97,10){\makebox(0,0)[l]{\LARGE{$\leftarrow$}}}
\put(80,0){\line(1,2){3}}
\put(92,2){\makebox(0,0)[c]{.......}}
\put(105,0){\line(1,2){3}}
\put(110,0){\circle*{1,3}}\put(110,0){\makebox(0,-6)[c]{\small $\bar j$}}
\put(115,0){\line(1,2){3}}
\put(125,2){\makebox(0,0)[c]{.......}}
\put(135,0){\line(1,2){3}}
\put(140,0){\circle*{1,3}}\put(140,0){\makebox(0,-6)[c]{\small $j$}}
\qbezier(110,0)(125,10)(140,0)
\put(95,-4){\makebox(0,0)[c]{\tiny $l_j-\gamma_j$ nestings}}
\put(125,-4){\makebox(1,0)[c]{\tiny $\gamma_j-1$ crossings}}
\put(30,-10){\makebox(0,0)[c]{\small $D_j$}}\put(110,-10){\makebox(0,0)[c]{\small $D_j'$}}
\end{picture}
}
\end{center}
\caption{Counting of $\cro(\pi;j)$ and $\cro(\varphi(\pi);j)$}
\end{figure}
Suppose that
  $j$ is linked  with
 the $\gamma_j$-th
vacant vertex $\bar j$ of $D_{j-1}$ in $D_j$
(resp. $D_{j-1}'(\pi)$ in $D_j'$).
Recall that the rank of vacant vertices is counted
from left to right in $D_{j-1}$ and  from right to left in $D_{j-1}'$.
\begin{itemize}
       \item Any  vacant vertex $\alpha$ on the left  of the vertex $z$
       in $D_j$ (resp. $D_j'$)
         will be linked to a vertex $\beta$ on the right  of the vertex $j$;
        thus  $(\alpha,\beta)$  will form a  nesting with $(\bar j,j)$  as an interior edge.
        Conversely,  if $(a,b)$ forms a  nesting  with  interior
        edge $(\bar j,j)$, then $a$ must be  a vacant vertex on the left  of the vertex $\bar j$
         in $D_{j}$ (resp. $D_j'$) .
         We deduce that $\ne(\pi;j)=\gamma_j-1$ and $\ne(\varphi(\pi);j)=l_j-\gamma_j$.
       \item Any  vacant vertex $\alpha$ on the right of the vertex $\bar j$ in $D_j$ (resp. $D_j'$)
        will be linked to a vertex $\beta$ on the right of the vertex $j$; thus  $(\alpha,\beta)$ will form
        a crossing  with  initial edge $(\bar j,j)$. Conversely, if $(a,b)$ forms  a crossing
        with  initial edge  $(\bar j,j)$ ,
        then the vertex $a$ must be  a vacant vertex on the right of the vertex $\bar j$ in $D_{j}$ (resp. $D_j'$).
        We deduce that $\cro(\pi;j)=l_j-\gamma_j$ and $\cro(\varphi(\pi);j)=\gamma_j-1$.
   \end{itemize}
 The proof is completed by comparing the above counting results.\qed

\section{Factorization of $\varphi$ via Charlier diagrams }

\subsection{Charlier diagrams}
 A {\it path} of length $n$ is a finite sequence
 $w=(s_{0},s_{1},\cdots,s_{n})$ of points
$s_{i}=(x_{i},y_{i})$ in the plan $\mathbb{Z}\times\mathbb{Z}$. A
step $(s_{i},s_{i+1})$ of $w$ is \emph{East} (resp.
\emph{North-East}  and \emph{South-East}) if
$s_{i+1}=(x_{i}+1,y_{i})$ (resp. $s_{i+1}=(x_{i}+1,y_{i}+1)$ and
$s_{i+1}=(x_{i}+1,y_{i}-1)$).
 The number $y_{i}$ is the {\emph height} of the step $(s_{i},s_{i+1})$.
 The integer $i+1$ is the {\emph index} of the step $(s_{i},s_{i+1})$.

  A {\it Motzkin path} is a path $w=(s_{0},s_{1},\cdots,s_{n})$
such that:
$s_{0}=(0,0)$ and $s_{n}=(n,0)$, each  step is  East or  North-East or South-East and
 $y_{i}\geq0$ for each $i$.
A {\it
bicolored Motzkin } (or BM) path  is a Motzkin path whose East
steps are
 colored with  \emph{red} or \emph{blue}.
A {\it restricted bicolored Motzkin} (or RBM) path  is a BM path
whose blue East steps are of height $>0$.

In the following, we shall write $BE$, $RE$, $NE$ and $SE$
as abbreviations of Blue East, Red East, North-East and South-East.

The type of $w$ is the 4-tuple $\lambda(w)=(\O(w), \C(w), \S(w),
\T(w))$, where  $\O(w)$ (resp. $\C(w)$, $\S(w)$, $\T(w)$) is the
set of indices of NE (resp. SE, RE, BE) steps of $w$.
 For instance, if $w$ is the path in Figure~7, then
 $$
 \lambda(w)=( \{1,2,5\} \,, \{7,10,11\} \,, \{4,8\} \,, \{3,6,9\}).
 $$
Denote by $M_b(n)$ (resp. $M_{rb}(n)$) the set of BM (resp. RBM) paths
of length $n$.

\begin{defn}
A \textit{Charlier diagram} of length $n$ is a  pair $h=(w,\xi)$
where $w=(s_{0},\ldots,s_{n})$ is a  RBM path and
$\xi=(\xi_{1},\ldots,\xi_{n})$ is a sequence of integers
such that $\xi_i=1$ if the $i$-th step is NE or RE,
and  $1\leq \xi_{i}\leq k$ if the $i$-th step is SE or
BE of height $k$.
\end{defn}

Let $\Gamma_n$ be the set of Charlier diagrams of length $n$. A Charlier diagram is given in Figure~7.
\begin{figure}[h]
{\setlength{\unitlength}{0.75mm}
\begin{center}
\begin{picture}(100,43)(0,-5)
\put(0,0){\line(1,0){110}}
%%%%%%%%%%%%%%%%%%%%%%%%%%%%
\put(0,0){\line(1,1){10}} \put(10,10){\line(1,1){10}}\put(20,20){\line(1,0){10}}\put(30,20){\line(1,0){10}}
\put(40,20){\line(1,1){10}}\put(50,30){\line(1,0){10}}\put(60,30){\line(1,-1){10}}\put(70,20){\line(1,0){10}}
\put(80,20){\line(1,0){10}}\put(90,20){\line(1,-1){10}}\put(100,10){\line(1,-1){10}}
%%%%%%%%%%%%%%%%%%%%%%%%%%%%%%%%%%%%%%%%%%%%%%%%%%%%%%%%%
\put(0,0){\circle*{2}}\put(10,10){\circle*{2}}\put(20,20){\circle*{2}}\put(30,20){\circle*{2}}\put(40,20){\circle*{2}}\put(50,30){\circle*{2}}
\put(60,30){\circle*{2}}\put(70,20){\circle*{2}}\put(80,20){\circle*{2}}\put(90,20){\circle*{2}}\put(100,10){\circle*{2}}\put(100,10){\circle*{2}}
%%%%%%%%%%%%%%%%%%%%%%%%%%%%%%%%%%%%%%%%%%%%%%%%%%%%%%%%%
\put(26,20){\makebox(0,3)[c]{\tiny blue}}\put(36,20){\makebox(0,3)[c]{\tiny red}}\put(56,30){\makebox(0,3)[c]{\tiny blue}}
\put(76,20){\makebox(0,3)[c]{\tiny red}}\put(86,20){\makebox(0,3)[c]{\tiny blue}}
%%%%%%%%%%%%%%%%%%%%%%%%%%%%%%%%%%%%%%%%%%%%%%%%%%%%%%%%%
\put(5,-5){\makebox(0,0)[c]{\tiny 1}}\put(15,-5){\makebox(0,0)[c]{\tiny 2}}\put(25,-5){\makebox(0,0)[c]{\tiny 3}}\put(35,-5){\makebox(0,0)[c]{\tiny 4}}
\put(45,-5){\makebox(0,0)[c]{\tiny 5}}\put(55,-5){\makebox(0,0)[c]{\tiny 6}}\put(65,-5){\makebox(0,0)[c]{\tiny 7}}\put(75,-5){\makebox(0,0)[c]{\tiny 8}}
\put(85,-5){\makebox(0,0)[c]{\tiny 9}}\put(95,-5){\makebox(0,0)[c]{\tiny 10}}\put(105,-5){\makebox(0,0)[c]{\tiny 11}}
%%%%%%%%%%%%%%%%%%%%%%%%%%%%%%%%%%%%%%%%%%%%%%%%%%%%%%%%%%%%
\put(-5,-5){\makebox(0,0)[c]{\tiny $i$}}\put(-5,-10){\makebox(0,0)[c]{\tiny $\xi_{i}$}}
\put(25,-10){\makebox(0,0)[c]{\tiny 2}}\put(55,-10){\makebox(0,0)[c]{\tiny 3}}\put(65,-10){\makebox(0,0)[c]{\tiny 2}}
\put(85,-10){\makebox(0,0)[c]{\tiny 1}}\put(95,-10){\makebox(0,0)[c]{\tiny 2}}\put(105,-10){\makebox(0,0)[c]{\tiny 1}}
\end{picture}
\end{center}}
\caption{\small a Charlier diagram of length 11}
\end{figure}
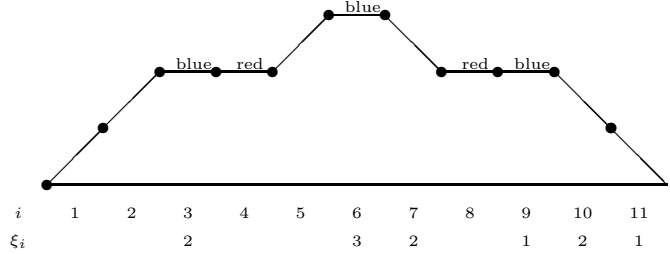

There is a well-known bijection (see \cite{Fl, Vi}) from
$\Gamma_n$ to $\Pi_n$. For our purpose, we present two variants
$\varphi_{l}$ and $\varphi_{r}$ of this bijection, which keep the
track of crossings and nestings.

Let $\pi$ be a partition. We denote
respectively by $\sg(\pi)$, $\bl(\pi)$ and $\tran(\pi)$ the
numbers of  singletons, blocks whose cardinal is $\geq2$ and
transients of the partition $\pi$.

 A partition $\pi$ of $[n]$ is
completely determined by its type $\lambda=(\O,\C,\S,\T)$ and the
integers $\gamma_i(\pi)$, $i\in \C\cup\T$. The description of the
bijection $\varphi_l$ is based on this fact. Given a Charlier
diagram $h=(w,\xi)$ of length $n$, we define the partition
$\pi=\varphi_l(h)$ as follows: the type of $\pi$ is that of $w$
and $\gamma_j(\pi):=\xi_j$, for $j\in \C\cup\T$.

In the definition of $\varphi_l$,
we take $\xi_j$ as the rank from left to right of the  vacant vertex linked
to $j$ in the $j$-th step of the construction of $\pi$. If we
we take $\xi_j$ as the rank from  right to left of the  vacant vertex linked
to $j$ in the $j$-th step of the construction of $\pi$, then we get the bijection $\varphi_r$.
That is,  we have $\varphi_r:=\varphi\,\circ\,\varphi_l$.
 In other words, the following diagram is commutative.
  \begin{figure}[h]
{\setlength{\unitlength}{0.60 mm}
\begin{center}
\begin{picture}(65,35)(0,-25)
\put(5,0){\line(1,0){60}}\put(5,0){\vector(1,0){60}}
\put(5,-3){\line(1,-1){25}}\put(5,-3){\vector(1,-1){25}}
\put(40,-28){\line(1,1){25}}\put(40,-28){\vector(1,1){25}}
\put(-1,0){\makebox(0,0)[c]{$\Pi_n$}}\put(71,0){\makebox(3,0)[c]{$\Pi_n$}}
\put(35,-30){\makebox(0,0)[c]{$\Gamma_n$}}\put(35,4){\makebox(0,0)[c]{$\varphi$}}
\put(18,-15){\makebox(-20,0)[c]{$\varphi_r^{-1}$}}\put(57,-15){\makebox(7,0)[c]{$\varphi_l$}}
\end{picture}
\end{center}}
\caption{factorization of $\varphi$}
\end{figure}
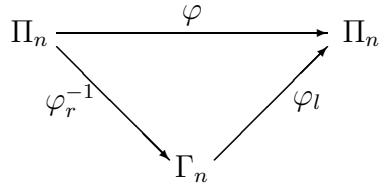
\goodbreak
The following result is clear (cf. \cite{Fl,Vi}).
\begin{prop}
The mapping $\varphi_l$ (resp. $\varphi_r):\Gamma_n\to\Pi_n$ is a bijection.
Moreover, if $h=(w,\xi)$ is a Charlier diagram and $\pi=\varphi_{r}(h)$ (or
$\pi=\varphi_{l}(h)$), then $\sg(\pi)$ (resp. $\bl(\pi)$ and
$\tran(\pi)$)
 is equal to the number of red East (resp. North-East and blue East ) steps of
 $w$.
\end{prop}

For instance, if $h=(w,\xi)$ is the Charlier diagram of
Figure~7, the construction of $\varphi_{l}(h)$ (resp.
$\varphi_{r}(h)$) correspond with the traces sequence $D_i(\pi)$
(resp. $D'_i$) in Figure~5. So,
$\varphi_l(h)=\{1,9,10\}-\{2,3,7\}-\{4\}-\{5,6,11\}-\{8\}$ and
$\varphi_r(h)=\{1,3,10\}-\{2,6,9,11\}-\{4\}-\{5,7\}-\{8\}$.

\begin{prop}\label{prop:psi}
 Let $h=(w,\xi)$ be a Charlier diagram such that the $j$-th step
of $w$ is blue East or South-East of height $k$, then
\begin{align*}
\cro(\varphi_{r}(h);j)&=\ne(\varphi_{l}(h);j)=\xi_{j}-1\\
\ne(\varphi_{r}(h);j)&=\cro(\varphi_{l}(h);j)=k-\xi_{j}
\end{align*}
\end{prop}

\pf This follows  from the  proof of Lemma~\ref{lem:principal} by
replacing $\varphi_{l}(h)$ by $\pi$, $\varphi_{r}(h)$ by
$\varphi(\pi)$, $l_j$ by $k$ and $\gamma_j$ by $\xi_j$. \qed

 A  partition $\pi$ is \emph{noncrossing} (resp. \emph{nonnesting})
if $\cro(\pi)=0$ (resp.  $\ne(\pi)=0$). Let $NC_n$ (resp. $NN_n$) be the set of noncrossing (resp. nonnesting)
partitions of $[n]$.
\begin{cor}\label{cor:noncr}
Let \textbf{1} denote the $n$-tuple $(1,1,\ldots,1)$. Then
\begin{itemize}
\item[(i)] The mapping  $w \mapsto \varphi_r ((w,\textbf{1}))$  is a bijection from $M_{rb}(n)$ to $NC_n$.
 \item[(ii)] The mapping  $w \mapsto \varphi_l ((w,\textbf{1}))$ is a bijection from $M_{rb}(n)$ to $NN_n$.
 \end{itemize}
 \end{cor}

\pf Let $h=(w,\xi)$ a restricted diagram and suppose that the $j$-th
step of $w$ is blue East or South-East. Then,
Proposition~\ref{prop:psi} implies that
$\cro(\varphi_{r}(h);j)=\ne(\varphi_{l}(h);j)=\xi_j-1$. Thus the partition
 $\varphi_{r}(h)$ (resp. $\varphi_{l}(h)$) is
noncrossing (resp. nonnesting) if and only if $\xi_i=1$
for each $i$.\qed

\begin{rmk}
Corollary~\ref {cor:noncr} gives another proof of the well-known fact (see \cite{Si}
and \cite[p.226]{Stan})
  that the cardinals
 of $NC_n$ and $NN_n$ are equal to the $n$-th Catalan number $C_n=\frac{1}{n+1}{2n\choose n}$.
Moreover, the mapping $\varphi=\varphi_l \circ \varphi_r^{-1}: NC_n\to NN_n$ is a bijection.
\end{rmk}

\subsection{Continued fraction expansions}
Consider the enumerating polynomial of $\Pi_n$:
 $$
 B_n(p,q,u_1,u_2,v)=\sum_{\pi\in\Pi_n}p^{\cro(\pi)}q^{\ne(\pi)}u_1^{sg(\pi)}u_2^{\bl(\pi)}v^{\tran(\pi)},
 $$
which  is a generalization of \emph{$n$-th Bell numbers}. Let
$$
[n]_{p,q}=\frac{p^n-q^n}{p-q},\qquad [n]_{q}=\frac{1-q^n}{1-q}.
$$
It follows from Proposition~\ref{prop:psi} that
\begin{align}\label{eq:key}
&B_n(p,q,u_1,u_2,v)\nonumber\\
&=\sum_{(w,\xi)\in \Gamma_n}\left(\prod_{j\in \O(w)}u_2\right)\left(\prod_{j\in \S(w)}u_1\right)
 \left(\prod_{j\in \C(w)}p^{\xi_j-1}q^{k_j-\xi_j}\right)
 \left(\prod_{j\in \T(w)}p^{\xi_j-1}q^{k_j-\xi_j}v\right),
  \end{align}
where $k_j$ is the height of the $j$-th step of $w$.

For any  BM path $w$,  define the \emph{weight} of a step of $w$ at
height $k$ by $u$ (resp. $[k]_{p,q}$, $d[k]_{p,q}(1-\delta_{0k})$,
$v$) if it is NE (resp. SE, BE, RE)
and  the weight $P(w)$ of $w$ as the product of weights of its steps.
We can rewrite the double sums in \eqref{eq:key} as a single sum on bicolored Motzkin paths:
$$
B_n(p,q,u_1,u_2,v)=\sum_{w \in M_b(n)}P(w).
$$
Applying a well-known result of Flajolet~\cite[Propositions 7A and 7B]{Fl}, we derive immediately the continued fraction expansion
from the above correspondence.
\begin{prop}\label{propCF}
The generating function
$\sum_{n\geq 0}B_n(p,q,u_1,u_2,v)z^n$ has the following continued fraction expansion:
\begin{align*}
\frac{1}{1-u_1z-\displaystyle\frac{u_2z^2}{1-(u_1+v)z-\displaystyle\frac{u_2[2]_{p,q}z^2}
{1-(u_1+[2]_{p,q}v)z-\displaystyle\frac{u_2[3]_{p,q}z^2}{1-(u_1+[3]_{p,q}v)z-
\displaystyle\frac{u_2[4]_{p,q}z^2}{\cdots}}}}}.
\end{align*}
\end{prop}
Note that
the $q=v=1$ case of Proposition~\ref{propCF} has been given by Biane~\cite{Bi}.
Taking $u_1=u_2=v=1$, we have:
\begin{cor}
 The generating function
$$\sum_{n\geq0}(\sum_{\pi\in\Pi_n}p^{\cro(\pi)}q^{ne(\pi)})z^n
=\sum_{n\geq0}(\sum_{\pi\in\Pi_n}q^{\cro(\pi)}p^{ne(\pi)})z^n$$ has the following continued fraction expansion:
\begin{align*}
\frac{1}{1-z-\displaystyle\frac{z^2}{1-([1]_{p,q}+1)z-\displaystyle\frac{[2]_{p,q}z^2}
{1-([2]_{p,q}+1)z-\displaystyle\frac{[3]_{p,q}z^2}{1-([3]_{p,q}+1)z-
\displaystyle\frac{[4]_{p,q}z^2}{\cdots}}}}}.
\end{align*}
\end{cor}

For any $\pi\in \Pi_n$, denote by $\ed(\pi)$ the number of edges of $\pi$.
Clearly we have  $\ed(\pi)=\bl(\pi)+\tran(\pi)$ and
$\cro(\pi)+\ne(\pi)+\al(\pi)={\ed(\pi)\choose 2}$.
Let
$$E_n(v,q):= \sum_{\pi\in\Pi_n}q^{\cro(\pi)+\ne(\pi)}v^{ed(\pi)}.
$$
Setting $p=q$, $u_1=1$ and $u_2=v$ in Proposition~\ref{propCF},  we get
\begin{cor}  The generating function $\sum_{n\geq0} E_n(v,q)z^n$ has the following continued fraction expansion:
\begin{align*}
\frac{1}{1-z-\displaystyle\frac{vz^2}{1-(1+v)z-\displaystyle\frac{2qvz^2}
{1-(2qv+1)z-\displaystyle\frac{3q^2vz^2}{1-(3q^2v+1)z-\displaystyle
\frac{4q^3vz^2}{\cdots}}}}}.
\end{align*}
\end{cor}

Let $E_n(v,q)=\sum_{k\geq 0}e_k(q)v^k$. Then
 $$
 F_n(q):= \sum_{\pi\in\Pi_n}q^{\al(\pi)}=\sum_{k\geq 0}q^{k\choose2} e_k(1/q).
 $$

Finally consider the enumerating polynomials of crossings and nestings of $\M_{2n}$:
$$L_n(p,q)=
\sum_{\alpha\in \mathcal M_{2n}}p^{\cro(\alpha)}q^{ne(\alpha)}
=\sum_{\alpha\in \mathcal M_{2n}}p^{ne(\alpha)}q^{\cro(\alpha)}.
$$
Setting $u_2=1$, $u_1=v=0$ in Proposition~\ref{propCF} and replacing $z^2$ by $z$ we get
\begin{prop}\label{prop:ptouchard}
$$
\sum_{n\geq0}L_n(p,q)z^n=
\frac{1}{1-\displaystyle\frac{z}{1-\displaystyle\frac{[2]_{p,q}z}
{1-\displaystyle\frac{[3]_{p,q}z}{1-\displaystyle\frac{[4]_{p,q}z}{\cdots}}}}}.
$$
\end{prop}
Note that the $p=1$ case of Proposition~\ref{prop:ptouchard} corresponds to a result of Touchard~\cite{Tou}.
Since
a matching of $[2n]$ has exactly $n$ edges, we get
$\cro(\alpha)+\ne(\alpha)+\al(\alpha)={n\choose 2}$ for any $\alpha\in \mathcal M_{2n}$.
Therefore
$$
T_n(q):=\sum_{\alpha\in \mathcal M_{2n}}q^{\al(\alpha)}=q^{{n\choose 2}}L_n(1/q,1/q).
$$
The first terms of the above sequences are given as follows:
$$
\begin{array}{ll}
   T_0(q)=T_1(q)=1    &                  L_0(p,q)=L_1(p,q)=1\\
   T_2(q)=2+q   &                 L_2(p,q)=1+p+q\\
   T_3(q)=6+4q+4q^2+q^3 &         L_3(p,q)=1+2p+2q+2pq+p^2+q^2+2p^2q+2pq^2+p^3+q^3.
\end{array}
$$

\renewcommand{\baselinestretch}{1}

\end{document}